\documentclass{svmult}
\usepackage{graphicx}

\begin{document}
\title*{Multiple testing with optimal individual tests \\ in Gaussian graphical model selection}
\author{Alexander P. Koldanov, Petr A. Koldanov, Panos Pardalos}
\institute{Alexander P. Koldanov \at National Research University
Higher School of Economics, Bolshaya Pecherskaya 25/12, Nizhny
Novgorod, 603155 Russia \email{akoldanov@hse.ru } \and Petr A.
Koldanov \at National Research University Higher School of
Economics, Bolshaya Pecherskaya 25/12, Nizhny Novgorod, 603155
Russia \email{pkoldanov@hse.ru} \and Panos Pardalos \at{University of Florida, Center for Applied Optimization,  Gainesville, FL, USA  }}
\maketitle
\abstract
{
Gaussian Graphical Model selection problem is considered. Concentration graph is identified by multiple decision procedure based on individual tests. 
Optimal unbiased individual tests are constructed. It is shown that optimal tests are equivalent to sample partial correlation tests. Associated multiple decision  procedure
is compared with standard procedure.
}
\keywords{Gaussian graphical model;  Bonferroni procedure;  Unbiased statistical tests; Exponential families; Tests of the Neyman structure.}






\section{Introduction}
Gaussian graphical model (GGM) is a useful and popular tool in biology and genetics \cite{bib1}, \cite{bib2}. This model is based on analysis of conditional independences between random variables which is equivalent to specified zeros among the set of partial correlation coefficients.  Some statistical procedures for Gaussian graphical model selection are proposed and investigated in 
\cite{bib3}, \cite{bib2}. Statistical procedures  in \cite{bib2} does not control family wise error rate (FWER). To overcome this issue multiple single-step and stepwise statistical procedures were recently introduced (see \cite{bib4}, \cite{bib5}, \cite{bib6} and references therein). Majority of multiple testing  procedures are based on statistics for individual hypotheses testing. Usually, sample partial correlations are used as natural test statistics for testing individual hypotheses. However, as far as we know, 
properties of optimality of these tests were not well investigated so far.  In the present paper we construct an optimal unbiased test of the Neyman structure for testing individual hypotheses in GGM selection and prove that these tests are equivalent to sample partial correlation tests. Numerical experiments are conducted to compare the multiple testing procedure with optimal individual tests with standard concentration graph identification 
procedure based on Fisher z-transform of sample partial correlations. 

The paper is organized as follows. Section \ref{Basic notations} contains basic definition and problem statement. In the Section \ref{Tests of the Neyman structure} a general description of the tests of the Neyman structure is given. In the Section \ref{Optimal individual test} optimal  tests for testing individual hypotheses are constructed. In the Section \ref{Optimality of the sample partial correlation tests} the main result of optimality of the sample partial correlation test is proved. Multiple testing procedure with optimal individual tests is presented in \ref{Multiple testing procedures} and it is compared with standard Bonferroni procedure. Concluding remarks are given in the Section \ref{conclusion}. 

\section{Basic notations and problem statement}\label{Basic notations}

Let random vector $X=(X_1,X_2,\ldots,X_N)$ has a multivariate Gaussian distribution $N(\mu, \Lambda)$, where $\mu=(\mu_1,\mu_2,\ldots, \mu_N)$ is the vector of means and 
$\Lambda=(\sigma_{i,j})$ is the covariance matrix, $\sigma_{i,j}=\mbox{cov}(X_i,X_j)$, $i,j=1,2,\ldots,N$. Denote by $\rho_{i,j}$ Pearson correlation between random variables $X_i$, $X_j$, 
$\rho_{i,j}=\sigma_{i,j}/\sqrt{\sigma_{i,i}\sigma_{j,j}}$. Let $x(t)$, $t=1,2,\ldots,n$  be a sample of the size $n$ from distribution $X$. Let
$$ 
s_{i,j}=\frac{1}{n} \Sigma_{t=1}^n(x_{i}(t)-\overline{x_{i}})(x_{j}(t)-\overline{x_{j}}), \ \ \overline{x_i}=(1/n)\sum_{t=1}^n x_i(t)
$$
be the sample covariance between $X_i$, $X_j$. Let $S=(s_{i,j})$ be the matrix of sample covariances. Denote by $r_{i,j}$ the sample Pearson correlation between random variables $X_i$, $X_j$, $r_{i,j}=s_{i,j}/\sqrt{s_{i,i}s_{j,j}}$. It is known \cite{bib7}, that the statistics $\overline{x}=(\overline{x_{1}},\overline{x_{2}},...,\overline{x_{N}})$ and $S=(s_{i,j})$ are sufficient for 
multivariate Gaussian distribution. The inverse matrix for $\Lambda$,  $\Lambda^{-1}=(\sigma^{i,j})$ is known as concentration or precision matrix for the distribution $X$.  

Consider the set $\cal{G}$ of all $N \times N$ symmetric matrices  $G=(g_{i,j})$ with $g_{i,j} \in \{0,1\}$, $i,j=1,2,\ldots,N$, $g_{i,i}=0$, $i=1,2,\ldots,N$. Matrices $G \in \cal{G}$ represent adjacency 
matrices of all simple undirected graphs with $N$ vertices. Total number of matrices in $\cal{G}$ equals to $L=2^M$ with $M=N(N-1)/2$. The problem of Gaussian graphical model
selection can be formulated now as a multiple decision problem of selecting one from a set of $L$ hypotheses:
\begin{equation}\label{N_hypotheses}
H_G: \sigma^{i,j}=0,\mbox{ if } g_{i,j}=0, \ \ \sigma^{i,j}\neq 0,\mbox{ if } g_{i,j}=1; \ \ i \neq j
\end{equation}
For the multiple decision problem (\ref{N_hypotheses}) we introduce the following set of individual hypotheses:
\begin{equation}\label{Individual hypotheses}
h_{i, j}:\sigma^{i,j}=0  \mbox{ vs } \
k_{i, j}:\sigma^{i,j}\neq 0, \\
\ i,\ j =1,2,\ldots, N, \ i \neq j.
\end{equation}
In this paper we consider individual tests for (\ref{Individual hypotheses}) of the following form
\begin{equation}\label{Individual tests}
\varphi_{i,j}(x)= \left \{
\begin{array}{lll}
0, & \mbox{if} & c^1_{i,j} < T_{i,j} < c^2_{i,j} \\
1, & \mbox{otherwise} & 
\end{array}
\right.
\end{equation}
where $T_{i,j}(x)$, are the tests statistics. There are many multiple testing procedures based on simultaneous use of individual tests statistics: Bonferroni, Holm, Hochberg, Simes, T-max procedures and others. 
In this paper we concentrate on Bonferroni multiple testing procedures with optimal individual tests. 

\section{Tests of the Neyman structure}\label{Tests of the Neyman structure}

To construct optimal individual tests for the problem (\ref{Individual hypotheses}) we use a tests of the Neyman structure \cite{bib8}, \cite{bib9}. 
Let $f(x;\theta)$ be the density of the  exponential family:
\begin{equation} \label{exp}
f(x;\theta)=c(\theta)exp(\sum_{j=1}^M\theta_jT_j(x))m(x)
\end{equation}
where $c(\theta)$ is a function defined in the parameters space,
$m(x)$, $T_j(x)$ are functions defined in the sample space $R^{N \times n}$, and
$T_j(X)$ are the sufficient statistics for $\theta_j,j=1,\ldots,M$.

Suppose that  individual  hypotheses  have the form:
\begin{equation}\label{generating hypotheses}
h_j:\theta_j=\theta^0_j \mbox{  vs  } k_j:\theta_j \neq\theta^0_j, \ \ j=1,2,\ldots,M
\end{equation}
where $\theta^0_j$ are  fixed. 

The optimal unbiased tests for individual hypotheses (\ref{generating hypotheses}) are (see \cite{bib8} ch. 4, theorem 4.4.1):
\begin{equation}\label{tsn}
\varphi_j= \left\{
\begin{array}{lc}
0, & \mbox{if} \ \ c^1_j(t_1,\ldots,t_{j-1},t_{j+1},\ldots,t_M)<t_j<c^2_j(t_1,\ldots,t_{j-1},t_{j+1},\ldots,t_M)\\
1, & \mbox{  otherwise  }
 \end{array}
\right.
\end{equation}
where $t_i=T_i(x),i=1,\ldots,M$ and constants $c^1_j$, $c^2_j$ are defined from the equations
\begin{equation}\label{neymstruc1}
\int_{c^1_j}^{c^2_j}f(t_j,\theta^0_j|T_i=t_i,i=1,\ldots,M;i \neq j)dt_j=1-\alpha_j
\end{equation}
and
\begin{equation}\label{neymstruc2}
\begin{array}{l}
\int_{-\infty}^{c^1_j}t_jf(t_j,\theta^0_j|T_i=t_i,i=1,\ldots,M;i\neq j)dt_j+\\
+\int_{c^2_j}^{+\infty}t_jf(t_j,\theta^0_j|T_i=t_i,i=1,\ldots,M;i\neq j)dt_j=\\
=\alpha_j\int_{-\infty}^{+\infty}t_jf(t_j,\theta^0_j|T_i=t_i,i=1,\ldots,M;i\neq
j)dt_j
\end{array}
\end{equation}
where $f(t_j,\theta^0_j|T_i=t_i,i=1,\ldots,M;i\neq j)$ is the
density of conditional distribution of statistics $T_j$  and $\alpha_{j}$ is the significance  level of
the test.

\section{Optimal individual tests}\label{Optimal individual test}
Now we construct the optimal test in the class of unbiased tests
for individual hypothesis (\ref{Individual hypotheses}). This construction is based on the tests of the Neyman structure. 
Joint distribution  of sufficient statistics $s_{k,l}$, $k,l = 1,2,\ldots,N$ is given by Wishart density function \cite{bib7}:
$$
f(\{s_{k,l}\})=\displaystyle \frac{ [\det (\sigma^{k,l})]^{n/2}
\times [\det(s_{k,l})]^{(n-N-2)/2}\times \exp[-(1/2)\sum_k \sum_l
s_{k,l} \sigma^{k,l}]} {2^{(Nn/2)}\times \pi^{N(N-1)/4} \times
\Gamma(n/2)\Gamma((n-1)/2)\cdots\Gamma((n-N+1)/2)}
$$
if the  matrix $(s_{k,l})$ is positive definite, and $f(\{s_{k,l}\})=0$ otherwise.  For a fixed $(i,j)$ one has:
$$
f(\{s_{k,l}\})=\displaystyle C(\{\sigma^{k,l}\}) 
\exp[-\sigma^{i,j}s_{i,j} - \frac{1}{2} \sum_{(k,l)\neq
(i,j);(k,l)\neq(j,i)} s_{k,l} \sigma^{k,l}]  m(\{s_{k,l}\})
$$
where
$$
C(\{\sigma^{k,l}\})=c_1^{-1}[\det (\sigma^{k,l})]^{n/2}
$$
$$
c_1=2^{(Nn/2)}\times \pi^{N(N-1)/4}\times
\Gamma(n/2)\Gamma((n-1)/2)\cdots\Gamma((n-N+1)/2)
$$
$$ m(\{s_{k,l}\})=[\det(s_{k,l})]^{(n-N-2)/2}
$$

According to (\ref{tsn}) the optimal test in the class of unbiased
tests for the individual hypothesis (\ref{Individual hypotheses}) has the form:
\begin{equation}\label{Nstructure}
\varphi_{i, j}(\{s_{kl}\})=\left\{\begin{array}{rl}
 \ 0, &\mbox{}\: if \:  c_{i,j}^1(\{s_{k,l}\})<s_{ij}<c_{i,j}^2 (\{s_{k,l}\}),\  (k,l)\neq (i,j)\\
 \ 1, &\mbox{}\: if \: s_{ij}\leq c_{i,j}^1(\{s_{kl}\})\mbox{ or } s_{i,j}\geq c_{i,j}^2(\{s_{k,l}\}),\  (k,l)\neq (i,j)
 \end{array}\right.
\end{equation}
where according to (\ref{neymstruc1}),(\ref{neymstruc2}) the critical values $c_{i,j}$  are defined from the equations
\begin{equation}\label{threshold1}
\displaystyle \frac{\int_{I \cap [c_{i,j}^1;c_{i,j}^2]}
 [\det(s_{k,l})]^{(n-N-2)/2}  ds_{i,j}}
{\int_{I}  [\det(s_{k,l})]^{(n-N-2)/2}
ds_{i,j}} =1-\alpha_{i,j}
\end{equation}
\begin{equation}\label{threshold2}
\begin{array}{l}
\displaystyle \int_{I \cap [-\infty;c_{i,j}^1]}
s_{i,j}[\det(s_{k,l})]^{(n-N-2)/2}
ds_{i,j}+\\
+\displaystyle \int_{I \cap [c_{i,j}^2;+\infty]}
s_{i,j} [\det(s_{k,l})]^{(n-N-2)/2}
ds_{i,j}=\\
 =\alpha_{i,j}\int_I s_{i,j}[\det(s_{k,l})]^{(n-N-2)/2} ds_{i,j}
\end{array}
\end{equation}
where $I$ is the interval of values of $s_{i,j}$ such that the
matrix $(s_{k,l})$ is positive definite and  $\alpha_{i,j}$ is the
 significance level of the tests.

Consider $\det(s_{k, l})$ as a function of the variable $s_{i,j}$. This determinant is a quadratic polynomial of $s_{i,j}$:
\begin{equation}\label{determinant_equation}
\det(s_{k,l})=-as_{i,j}^2+bs_{i,j}+c
\end{equation} 

Let $K=(n-N-2)/2$. Denote by $x_1,x_2$  the roots of the equation $-ax^2+bx+c=0$. One has with the change of variable $x=x_1+(x_2-x_1)u$:
$$
\int_f^d(ax^2-bx-c)^Kdx=\displaystyle (-1)^Ka^K(x_2-x_1)^{2K+1}\int_{\frac{f-x_1}{x_2-x_1}}^{\frac{d-x_1}{x_2-x_1}}u^K(1-u)^Kdu
$$
Therefore the equation (\ref{threshold1}) takes the form:
\begin{equation}\label{threshold_neyman_structure_1}
\displaystyle \int_{\frac{c^1-x_1}{x_2-x_1}}^{\frac{c^1-x_1}{x_2-x_1}}u^K(1-u)^Kdu=(1-\alpha)\int_0^1 u^K(1-u)^K du
\end{equation}
or
\begin{equation}\label{beta_function}
\displaystyle \frac{\Gamma(2K+2)}{\Gamma(K+1)\Gamma(K+1)}\int_{\frac{c^1-x_1}{x_2-x_1}}^{\frac{c^1-x_1}{x_2-x_1}}u^K(1-u)^Kdu=(1-\alpha)
\end{equation}
It means that conditional distribution of $s_{i,j}$ when all other $s_{k,l}$ are fixed is the beta distribution $Be(K+1,K+1)$.

Beta distribution $Be(K+1,K+1)$ is symmetric with respect to the point $\frac{1}{2}$. Therefore the significance level condition (\ref{threshold1}) 
and unbiasedness condition (\ref{threshold2}) are satisfied if and only if:
$$
\displaystyle \frac{c^2-x_1}{x_2-x_1}=1 - \frac{c^1-x_1}{x_2-x_1}
$$  

Let $q$ be the $\frac{\alpha}{2}(\alpha<\frac{1}{2})$-quantile of beta distribution $Be(K+1,K+1)$, i.e. $F_{Be}(q)=\frac{\alpha}{2}$. Then thresholds  $c^1$, $c^2$ are defined by: 
\begin{equation}\label{thresholds_from_beta_distr}
\begin{array}{c}
c^1=x_1+(x_2-x_1)q \\
c^2=x_2-(x_2-x_1)q
\end{array}
\end{equation}

Finally, the test of the Neyman structure accepts the hypothesis $\sigma^{i,j}=0$ vs alternative $\sigma^{i,j}\neq 0$ iff:
\begin{equation}\label{Neyman_structure_test}
c^1_{i,j} < s_{i,j} < c^2_{i,j} 
\end{equation} 
where 
$$
\begin{array}{l}
c^1_{i,j}=x_1(i,j)+(x_2(i,j)-x_1(i,j))q(\alpha_{i,j}) \\
c^1_{i,j}=x_2(i,j)-(x_2(i,j)-x_1(i,j))q(\alpha_{i,j})
\end{array}
$$
and $F_{Be}(q(\alpha_{i,j}))=\alpha_{i,j}/2$, $x_1,x_2$ are defined from the equation $ax^2-bx-c=0$ with $a,b,c$ defined in (\ref{determinant_equation}). 
The test (\ref{Neyman_structure_test}) can be written in the following form
\begin{equation}\label{Neyman_structure_q}
\varphi_{i,j}=\left\{\ 
\begin{array}{ll} 
0, & \displaystyle 2q-1 < \frac{as_{i,j}-\frac{b}{2}}{\sqrt{\frac{b^2}{4}+ac}} < 1-2q \\
1, & \mbox{otherwise}
\end{array}\right.
\end{equation}
where $a=a_{i,j},b=b_{i,j},c=c_{i,j}$ are defined in (\ref{determinant_equation}). 

\section{Optimality of the sample partial correlation tests}\label{Optimality of the sample partial correlation tests}

It is known \cite{bib1} that hypothesis $\sigma^{i,j}=0$ is equivalent to the hypothesis  $\rho^{i,j}=0$, 
where $\rho^{i,j}$ is the saturated partial correlation between $X_i$ and $X_j$:
$$
\rho^{i,j}=\frac{- \Lambda_{i,j}}{\sqrt{\Lambda_{i,i}\Lambda_{j,j}}}
$$
where for a given matrix $A=(a_{k,l})$ we denote by $A_{i,j}$  the cofactor of the element $a_{i,j}$. 
Denote by $r^{i,j}$ a sample partial correlation 
$$
r^{i,j}=\frac{-S_{i,j}}{\sqrt{S_{i,i}S_{j,j}}}
$$
where $S_{i,j}$ is the cofactor of the element $s_{i,j}$ in the matrix of sample covariances $S$.

Test for testing hypothesis $\rho^{i,j}=0$ has the form \cite{bib7}:
\begin{equation}\label{Partial_correlation_test}
\varphi_{i,j}=\left\{\ 
\begin{array}{ll} 
0,&|r^{i,j}|\leq c_{i,j}\\
1,&|r^{i,j}|> c_{i,j}
\end{array}\right.
\end{equation} 
where 
$c_{i,j}$ is $(1-\alpha_{i,j}/2)$-quantile of the distribution with the following density function
$$
f(x)=\displaystyle \frac{1}{\sqrt{\pi}}\frac{\Gamma(n-N+1)/2)}{\Gamma((n-N)/2)}(1-x^2)^{(n-N-2)/2}, \ \ \ -1 \leq x \leq 1
$$
In practical applications  the following Fisher transformation is used:
$$
z_{i,j}=\frac{\sqrt{n}}{2} \ln\left(\frac{1+r^{i,j}}{1-r^{i,j}}\right)
$$
For the case $\rho^{i,j}=0$ statistic $z_{i,j}$ asymptotically has standard Gaussian distribution. 
That is why the following test is largely used \cite{bib7}, \cite{bib5}:
\begin{equation}\label{Fisher_test}
\varphi_{i,j}=\left\{\ 
\begin{array}{ll} 
0,&|z_{i,j}|\leq c_{i,j}\\
1,&|z_{i,j}|> c_{i,j}
\end{array}\right.
\end{equation} 
where the constant $c_{i,j}$ is $(1-\alpha_{i,j}/2)$-quantile  of standard Gaussian distribution. 

In this section we prove that optimal unbiased test  (\ref{Neyman_structure_q}) is equivalent to the test (\ref{Partial_correlation_test}).

\noindent
{\bf Theorem:} Sample partial correlation test (\ref{Partial_correlation_test}) is optimal in the class of unbiased tests for testing hypothesis 
$\rho^{i,j}=0$ vs $\rho^{i,j} \neq 0$.

\noindent
{\bf Proof:} it is sufficient to prove that
\begin{equation}\label{equality}
\displaystyle \frac{S_{i,j}}{\sqrt{S_{i,i}S_{j,j}}}=\frac{as_{i,j}-\frac{b}{2}}{\sqrt{\frac{b^2}{4}+ac}}
\end{equation}
To prove this equation we introduce some notations. Let $A=(a_{k,l})$ be an $(N \times N)$ symmetric matrix. Fix $i<j$, $i,j =1,2,\ldots,N$.  
Denote by $A(x)$ the matrix obtained from $A$ by replacing the elements $a_{i,j}$ and $a_{j,i}$  by $x$. 
Denote by $A_{i,j}(x)$ the cofactor of the element $(i,j)$ in the matrix $A(x)$. Then the following statement is true

\noindent
{\bf Lemma:} One has $[\mbox{det}A(x)]' = -2A_{i,j}(x)$.

\noindent
{\bf Proof of the Lemma:}  one has from the general Laplace decomposition of $\det A(x)$ by two rows $i$ and $j$:
$$
\det(A(x))=\det \left( 
\begin{array}{ll}
a_{i,i} & x \\
x & a_{j,j} \\
\end{array}  
\right) A_{\{i,j\},\{i,j\}} + \sum_{k<j, k \neq i} \det \left( 
\begin{array}{ll}
a_{i,k} & x \\
a_{j,k} & a_{j,j} \\
\end{array}  
\right)A_{\{i,j\},\{k,j\}}+
$$
$$ + \sum_{k>j} \det \left( 
\begin{array}{ll}
x & a_{i,k} \\
a_{j,j} & a_{j,k} \\
\end{array}  
\right)A_{\{i,j\},\{j,k\}}+ \sum_{k<i} \det \left( 
\begin{array}{ll}
a_{i,k} & a_{i,i} \\
a_{j,k} & x \\
\end{array}  
\right)A_{\{i,j\},\{k,i\}}+ 
$$
$$
\sum_{k>i, k \neq j} \det \left( 
\begin{array}{ll}
a_{i,i} & a_{i,k}\\
x & a_{j,k} \\
\end{array}  
\right)A_{\{i,j\},\{i,k\}}+\sum_{k<l, k,l \neq i,j} \det \left( 
\begin{array}{ll}
a_{i,k} & a_{i,l}\\
a_{j,k} & a_{j,l} \\
\end{array}  
\right)A_{\{i,j\},\{k,l\}}
$$
where $A_{\{i,j\},\{k,l\}}$ is the cofactor of the matrix 
$\left( 
\begin{array}{ll}
a_{i,k} & a_{i,l}\\
a_{j,k} & a_{j,l} \\
\end{array}  
\right)
$
in the matrix $A$. Taking the derivative of $\det A(x)$ one get
$$
[\det(A(x))]'=-2xA_{\{i,j\},\{i,j\}}-\sum_{k<j, k \neq i} a_{j,k}A_{\{i,j\},\{k,j\}}+
\sum_{k>j} a_{j,k}A_{\{i,j\},\{j,k\}}+
$$
$$
+\sum_{k<i} a_{k,i}A_{\{i,j\},\{k,i\}}-
\sum_{k>i, k \neq j} a_{k,i}A_{\{i,j\},\{k,i\}}=-2A_{i,j}(x)
$$
The last equation follows from the symmetry conditions $a_{k,l}=a_{l,k}$ and from Laplace decomposition of $A_{i,j}(x)$ by the row $j$ and the column $i$. 
Lemma is proved.

Now we come back to the proof of the theorem. One has $\det(S(x))=-ax^2+bx+c$, where $a,b,c$ are the same as in (\ref{determinant_equation}). Therefore by Lemma one has $[\det S(x)]'=-2ax+b=-2S_{i,j}(x)$, i.e. 
$S_{i,j}(x)=ax-b/2$. Let 
$x=s_{i,j}$ then $as_{i,j}-\frac{b}{2}=S_{i,j}$.  To prove the theorem it is sufficient to prove that $\sqrt{S_{i,i}S_{j,j}}=\sqrt{\frac{b^2}{4}+ac}$. 
Let $x_2=\frac{b+\sqrt{b^2+4ac}}{2a}$ be the maximum root of equation $ax^2-bx-c=0$. Then $ax_2-\frac{b}{2}=\sqrt{\frac{b^2}{4}+ac}$.  
Consider
$$
r^{i,j}(x)=\frac{-S_{i,j}(x)}{\sqrt{S_{i,i}S_{j,j}}}
$$
The value of $r^{i,j}(x)$ is a partial correlation associated with the covariance matrix $S(x)$. When $x$ is increasing from $x_1$ to $x_2$ then $r^{i,j}(x)$ is decreasing from $1$ to $-1$.
That is $r^{i,j}(x_2)=-1$, i.e.  $ax_2-\frac{b}{2}=\sqrt{S_{i,i}S_{j,j}}$. Therefore 
$$
\sqrt{S_{i,i}S_{j,j}}=\sqrt{\frac{b^2}{4}+ac}
$$
The Theorem is proved. Therefore the optimal unbiased test for testing hypothesis $\sigma^{i,j}=0$ vs $\sigma^{i,j} \neq 0$ can be written in the following form
\begin{equation}\label{Neyman_structure_final}
\varphi_{i,j}=\left\{\ 
\begin{array}{ll} 
0,& \displaystyle 2q-1 < r^{i,j} < 1-2q \\
1& \mbox{otherwise}
\end{array}\right.
\end{equation}
where $q$ is the $\frac{\alpha}{2}$-quantile of beta distribution $Be(K+1,K+1)$.

\section{Multiple testing procedures}\label{Multiple testing procedures}
Consider two multiple decision Bonferroni type statistical procedures for the problem (\ref{N_hypotheses}). 
First procedure will be based on standard individual tests (\ref{Fisher_test}). Second procedure will be based on the tests of the Neyman structure (\ref{Neyman_structure_q}).
Let $\Phi^{Fi}(x)$ be the matrix
\begin{equation}\label{test_for_N_hypotheses_Fisher}
\Phi^{Fi}(x)=\left(\begin{array}{cccc}
1,&\varphi^{Fi}_{1,2}(x),&\ldots,&\varphi^{Fi}_{1,N}(x)\\
\varphi^{Fi}_{2,1}(x),&1,&\ldots,&\varphi^{Fi}_{2,N}(x)\\
\ldots&\ldots&\ldots&\ldots\\
\varphi^{Fi}_{N,1}(x),&\varphi^{Fi}_{N,2}(x),&\ldots,&1\\
\end{array}\right).
\end{equation}
where $\varphi^{Fi}_{i,j}(x)$ are defined by (\ref{Fisher_test}),  and constants 
$c_{i,j}=c$ are $(1-\alpha/(N(N-1))$-quantiles of standard Gaussian distribution. In this case the probability of at least one error (FWER)
is asymptotically majorated by $\alpha$. 
Bonferroni type multiple decision statistical procedure based on standard individual tests (\ref{Fisher_test}) is given by  
\begin{equation}\label{mdp_another_form_Fisher}
\delta_1(x)=d_G, \  \mbox{iff} \  \Phi^{Fi}(x)=G
\end{equation}

\begin{figure}[!h]
\includegraphics[scale=.35]{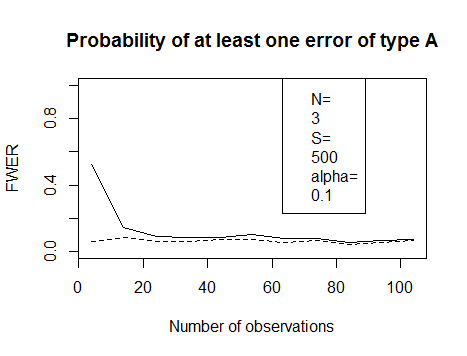}
\includegraphics[scale=.35]{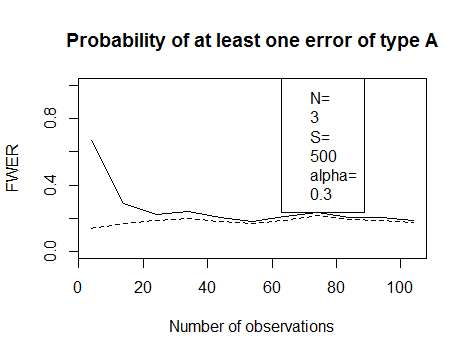}
\caption{FWER as a function of $n$. N=3. Left: significance level $\alpha=0.1$. Right: significance level $\alpha=0.3$ 
Solid line - procedure $\delta_1$ (standard). 
Dashed line - procedure $\delta_2$ (optimal). 
Horizontal axe represents the number of observations $n$.}
\label{N_3_diag}
\end{figure}

\begin{figure}[!h]
\includegraphics[scale=.35]{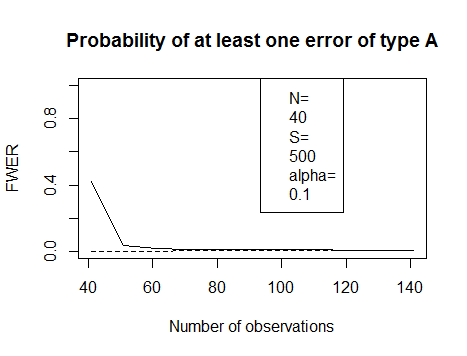}
\includegraphics[scale=.35]{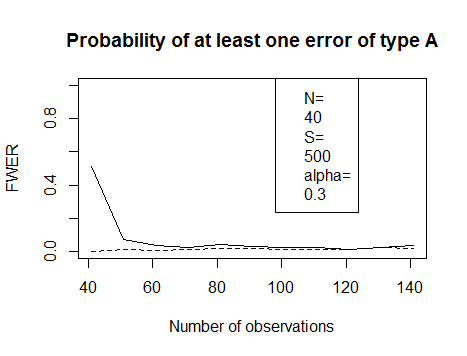}
\caption{FWER as a function of $n$. N=40. Left: significance level $\alpha=0.1$. Right: significance level $\alpha=0.3$ 
Solid line - procedure $\delta_1$ (standard). 
Dashed line - procedure $\delta_2$ (optimal). 
Horizontal axe represents the number of observations $n$.}
\label{N_40_diag}
\end{figure}

\begin{figure}[!h]
\includegraphics[scale=.35]{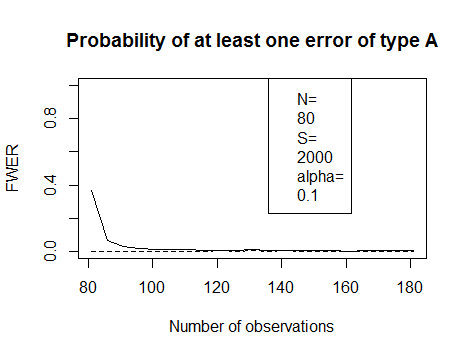}
\includegraphics[scale=.35]{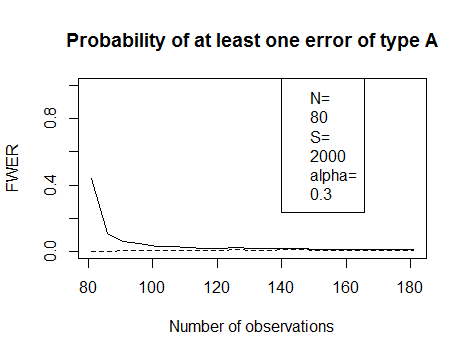}
\caption{FWER as a function of $n$. N=80. Left: significance level $\alpha=0.1$. Right: significance level $\alpha=0.3$ 
Solid line - procedure $\delta_1$ (standard). 
Dashed line - procedure $\delta_2$ (optimal). 
Horizontal axe represents the number of observations $n$.}
\label{N_80_diag}
\end{figure}

The second procedure $\delta_2$ is defined by the decision matrix $\Phi^{Ne}(x)$: 
\begin{equation}\label{test_for_N_hypotheses_Neyman_structure}
\Phi^{Ne}(x)=\left(\begin{array}{cccc}
1,&\varphi^{Ne}_{1,2}(x),&\ldots,&\varphi^{Ne}_{1,N}(x)\\
\varphi^{Ne}_{2,1}(x),&1,&\ldots,&\varphi^{Ne}_{2,N}(x)\\
\ldots&\ldots&\ldots&\ldots\\
\varphi^{Ne}_{N,1}(x),&\varphi^{Ne}_{N,2}(x),&\ldots,& 1\\
\end{array}\right).
\end{equation}
where $\varphi^{Ne}_{ij}(x)$ are defined by (\ref{Neyman_structure_final}) and 
$q$ is the $\frac{\alpha}{N(N-1)}$-quantile of beta distribution $Be(K+1,K+1)$. In this case the probability of at least one error (FWER)
is exactly majorated by $\alpha$. Two multiple decision procedures $\delta_1(x)$ and $\delta_2(x)$ have similar FWERs  when number of observations $n$ is sufficiently large. 
It is interesting to compare FWER for these procedures  for a small sample size $n$. The Figures \ref{N_3_diag}-\ref{N_80_diag} present a results of numerical experiments conducted for the case of 
diagonal matrix $\Lambda$. This is the case where one can expect a highest value of FWER. We use $S$ replications to estimate the probability of at least one false edge inclusion (FWER).  

One can note that standard procedure $\delta_1$ does not control the FWER for number of observation $n$ close to the matrix dimension $N$. 
It is interesting to note that  FWER for $\delta_1$ is rapidly decreasing starting from the value $n=N+1$,  it is stabilized from some value of $n=n(N)$ and
the difference $(n(N)-N)$ is decreasing in contrast with what one can expect.

\section{Concluding remarks}\label{conclusion}
Optimality of individual tests for individual hypotheses testing does not imply in general optimality of associated single-step multiple testing procedure. 
However if the losses from false decisions are supposed to be additive then in some cases it is possible to prove optimality of multiple testing procedure
\cite{bib10}, \cite{bib11}. Application of this approach for GGM selection will be a subject of forthcoming publication.

\begin{acknowledgement}
The work was conducted at National Research University Higher School of Economics, Laboratory of Algorithms and Technologies for Network Analysis. 
Partly supported by NRU HSE Scientific Fund and RFFI 14-01-00807.
\end{acknowledgement}




\end{document}